\documentclass{birkjour}
 \newtheorem{thm}{Theorem}[section]
 
 \newtheorem{lem}[thm]{Lemma}
 
 \theoremstyle{definition}
 
 \theoremstyle{remark}
 \newtheorem{rem}[thm]{Remark}
 
 \numberwithin{equation}{section}

\allowdisplaybreaks
\begin{document}

\title[Sz\'asz-P{\u{a}}lt{\u{a}}nea type operators involving Appell polynomials of class $A^2$]{Some Approximation Properties by Sz\'asz-P{\u{a}}lt{\u{a}}nea type Operators involving the Appell Polynomials of class $A^2$}

\author{Naokant Deo}

\address{
Department of Applied Mathematics,\\
Delhi Technological University,\\
Bawana Road Delhi-110042, INDIA}

\email{naokantdeo@dce.ac.in}


\author{Chandra Prakash\thanks{correspondingauthor}}
\address{Department of Applied Mathematics,\\
Delhi Technological University,\\
Bawana Road Delhi-110042, INDIA}
\email{chandu.math.bhu@gmail.com}

\author{D. K. Verma}
\address{Department of Mathematics,\\
Miranda House,\\
University of Delhi,\\
Delhi-110007, INDIA}
\email{durvesh.kv.du@gmail.com}
\subjclass{41A10; 41A25; 41A35; 41A36}

\keywords{Lipschitz type space, Voronvaskaja type asymptotic result, Modulus of continuity, Weighted approximation, Derivative of bounded Variation}

\date{December 26, 2022}
\begin{abstract}
This article contributes to the new summation of Sz\'asz operators with the help of Appell polynomials of class $A^{2}$. We verified Bohman-Korovkin's theorem and prove the convergence results like Lipschitz-type space, Voronvaskaja-type asymptotic formula, and modulus of continuity using the given operators. Furthermore, we have shown the weighted modulus of continuity and the derivative of bounded variation.
\end{abstract}

\maketitle

\section{Introduction}
One of the most important subfields in functional analysis is operator theory. Many continuous real-valued functions maintain their linearity on closed and constrained intervals. S. N. Bernstein provided a sequence of polynomials for the real-valued functions defined in 1912 on the finite closed interval [0, 1],
\begin{eqnarray*}
B_{n}(f;x)&=& \sum_{i=0}^{n}\binom{n}{i}x^{i}(1-x)^{n-i}f\left(\frac{i}{n}\right),\,\,\,\forall\,\, x\in [0,1],\,\,n\in \mathbb{N}.
\end{eqnarray*}
Using these polynomials, Bernstein also provided a simple proof for the Weierstrass approximation theorem. These finite summation polynomials were the focus of a large number of studies, which produced approximation findings for the finite intervals (\cite{CPDKND21,CPNDDK21,HGRP}). Sz\'asz compiled the Bernstein polynomials work on the infinite interval in 1950 and came up with the convergence findings. Moreover, these polynomials are significant in the area of approximation theory. The following are Sz\'asz polynomials:

\begin{eqnarray}\label{eq1.1}
S_{n}(f;x)&=& e^{-nx}\sum_{i=0}^{\infty}\frac{(nx)^{i}}{i!}f\left(\frac{i}{n}\right),
\end{eqnarray} where $f\in C[0,\infty)$, for each $x\in [0,\infty)$.
Jakimovski et al. [12] promoted Sz'asz operators with the use of the Appell polynomials in the context of the work done in the field of Sz'asz operators.
Consider an analytic function in the disc $|z|<\mathbb{R},\,(\mathbb{R}>1)$ and $g(z)=\displaystyle \sum_{i=0}^{\infty}a_{i}z^{i}\,\,(a_{0}\neq 0)$ with $g(1)\neq 0$. The Appell polynomials $p_{i}(x)$ contains generating function in the form $g(u)e^{ux}=\displaystyle \sum_{i=0}^{\infty}p_{i}(x)u^{i}$. Jakimovski and Leviatan build the following positive linear operators $p_{n}(f;x)$ defined as:
\begin{eqnarray}\label{eq1.2}
P_{n}(f;x)&=&\frac{e^{-nx}}{g(1)}\sum_{i=0}^{\infty}p_{i}(nx)f\left(\frac{i}{n}\right),\,\,\forall\,\, n\in \mathbb{N},
\end{eqnarray} and they deliberate the approximation results on the basis of Sz\'{a}sz operators (\ref{eq1.1}). In \cite{BW}, Wood  studied the approximation properties of the operators (\ref{eq1.2}) on $[0,\infty)$ iff $\frac{a_{i}}{g(1)}\geq 0$, for $i\in \mathbb{N}$. When $g(u)=1$, then operators (\ref{eq1.2}) reduces to (\ref{eq1.1}). Using the above certitude, Ciupa (\cite{AC1996,AC1994}) gave another variants of the operators and obtained the approximation properties by the goodness of Korovkin's theorem.
In $1974$, Ismail \cite{IMEH}, studied another generalization of Sz\'asz operators in the view of Sheffer polynomials and operators $P_{n}(f;x)$. Let $A(z)=\displaystyle \sum_{i=0}^{\infty}a_{i}z^{i},\,\,(a_{0}\neq 0)$ and $H(z)=\displaystyle \sum_{i=1}^{\infty}h_{i}z^{i}\,\,(h_{1}\neq 0)$ be analytic functions in the disc $|z|<\mathbb{R},\,\,(\mathbb{R}>1)$, where $a_{i}$ and $h_{i}$ in real. In view of the following assumptions
\begin{enumerate}
\item $p_{i}(x)\geq 0,$ for $x\in [0,\infty)$.
\item $A(1)\neq 0$ and, $H^{\prime}(1)=1$.
\end{enumerate} The set of polynomials $p_{i}(x);\,i\geq 0$ are sheffer polynomials $\Leftrightarrow$ the generating function of the form
\begin{eqnarray}\label{eq1.3}
A(\zeta)e^{xH(\zeta)}&=& \sum_{i=0}^{\infty}p_{i}(x)\zeta^{i},\,|\zeta|<\mathbb{R}.
\end{eqnarray} Finally, Ismail found the approximation properties of the given function
\begin{eqnarray}\label{eq1.4}
\Im_{n}(f;x)&=& \frac{e^{-nxH(1)}}{A(1)}\sum_{i=0}^{\infty}p_{i}(nx)f\left(\frac{i}{n}\right),\,n>0.
\end{eqnarray} Substitute $H(\zeta)=\zeta$ in the equation (\ref{eq1.4}), then the operators reduces to Leviatan operators (\ref{eq1.2}). And another way, we supposing $H({\zeta})=\zeta,\, A({\zeta})=1$, then the operators (\ref{eq1.4}) reduces to the Sz\'asz operators (\ref{eq1.1}).
Continuation on the above work, In $2015$, Sezgin Sucu and, Serhan Varma \cite{SSSV} obtained Stancu form of the $\Im_{n}$ operators (\ref{eq1.4}) and estimated their approximation properties. The Jakimovski-Leviatan-Stancu-Durrmeyer type operators were proposed by M. Mursaleen et al. In 2019 [17], they utilized the proposed operators to find the rth derivative of a function and to provide approximation qualities such as statistical convergence. Sz\'asz operators have been the subject of a great deal more study. Systems for computer algebra, data science, machine learning, and image processing all make substantial use of these operators.
Using the Appell polynomials of class $A2$, Varma and Sucu recently introduced the generalisation of Sz'asz operators as follows:
\begin{eqnarray}\label{eq1.5}
T_{n}(f;x)=\frac{1}{A(1)e^{nx}+B(1)e^{-nx}}\sum_{i=0}^{\infty}p_{i}(nx)f\left(\frac{i}{n}\right),
\end{eqnarray}with the restrictions $A(1)>0, B(1)\geq 0$ and $p_{i}(x)>0$ for all $i=0,1,\ldots$ The Appell polynomials $p_{i}(x)$ of class $A^{2}$ are obtained by the following generating function
\begin{eqnarray}\label{eq1.6}
A(\zeta)e^{x\zeta}+B(\zeta)e^{-x\zeta}=\sum_{i=0}^{\infty}p_{i}(x)\zeta^{i},
\end{eqnarray} where $A(\zeta)=\displaystyle \sum_{i=0}^{\infty}\frac{a_{i}}{i!}\zeta^{i}$ and $B(\zeta)=\displaystyle \sum_{i=0}^{\infty}\frac{b_{i}}{i!}\zeta^{i}$ are power series defined over the disc $|z|<\mathbb{R}\,\,(\mathbb{R}>1)$ with $a_{0}^2-b_{0}^2\neq 0$. Using recently proposed operators, they achieved approximation results by using the modulus of continuity, and they also examined at a specific sequence of operators that contained Gould-Hopper polynomials.

We studied several summation-integral type operators as a consequence of some research publications (\cite{CPNDDK, CPM77, NKGVG06,VGRN06}), and with the assistance of Varma and Sucu's work,  we present an innovative generalization of Sz\'asz-P{\u{a}}lt{\u{a}}nea operators utilizing the Appell polynomials of class $A^{2}$  on the interval $[0,\infty)$, define as
\begin{eqnarray}\label{eq1.7}
S_{n}^{\rho}(f;x)&=& \frac{1}{A(1)e^{nx}+B(1)e^{-nx}}\left[\sum_{i=1}^{\infty}p_{i}(nx)\int_{0}^{\infty}\phi_{n,i}^{\rho}(\zeta,c)f(\zeta)d\zeta\nonumber \right.\\
&&\left.+(a_{0}e^{nx}+b_{0}e^{-nx})f(0)\right]
\end{eqnarray}
\begin{eqnarray*}
\phi_{n,i}^\rho(\zeta,c)=\left\{
\begin{array}{c}
\dfrac{n\rho}{\Gamma{(i\rho)}}e^{-n\rho \zeta}(n\rho \zeta)^{i\rho-1}, c=0 \\
\dfrac{\Gamma\left(\frac{n\rho}{c}+i\rho\right)}{\Gamma(i\rho)\Gamma\left(\frac{n\rho}{c}\right)}\dfrac{c^{i\rho}\zeta^{i\rho-1}}{(1+c \zeta)^{\frac{n\rho}{c}+i\rho}},c=1,2,3,...
\end{array}
\right.\end{eqnarray*}
It can we easily observed by simple computation that:
\begin{eqnarray*}
\int_{0}^\infty\phi_{n,i}^\rho \left(\zeta,c\right)\zeta^r d\zeta=\left\{\begin{array}{ll}
\frac{\Gamma(i\rho+r)}{\Gamma(i\rho)}\dfrac{1}{\prod_{j=1}^r\left(n\rho-jc\right)},& r\neq0\\
1,& r=0.
\end{array}\right.
\end{eqnarray*}

On the other hand, given operators (\ref{eq1.7}) can be write as
\begin{eqnarray}\label{eq1.8}
S_{n}^{\rho}(f;x)&=& \int_{0}^{\infty}A_{n}(x,\zeta)f(\zeta)d\zeta,
\end{eqnarray} where

\begin{eqnarray*}
A_{n}(x,\zeta)=\frac{1}{A(1)e^{nx}+B(1)e^{-nx}}\left[(a_{0}e^{nx}+b_{0}e^{-nx})\delta(t)+\sum_{i=1}^{\infty}p_{i}(nx)\phi_{n,i}^{\rho}(\zeta,c)\right],
\end{eqnarray*} and $\delta(t)$ is the Dirac-delta function.\\
For $\beta>0$, let $C_{\beta}[0,\infty)=\{f\in C[0,\infty):\,f(\zeta)\leq N(1+\zeta^{\beta})\}$, for some $M>0$ endowed with the norm $$\|f\|_{\beta}=\sup_{\zeta\in [0,\infty)}\frac{|f(\zeta)|}{(1+\zeta^{\beta})}.$$
\begin{rem}
\begin{enumerate}
\item[(a)] For $A(\zeta)=1$ and $B(\zeta)=0$ in the equations (\ref{eq1.5}) and (\ref{eq1.7}), we obtain the famous Sz\'asz operators \cite{OZ} and Sz\'asz-P{\u{a}}lt{\u{a}}nea operators \cite{RP}.
\item[(b)]  For $A(\zeta)=1, \,B(\zeta)=0$ and $\rho=1$, the operators \eqref{eq1.7} reduces to Philips operators \cite{PHLP}.
\end{enumerate}
\end{rem}

\section{Lemmas}
In the sequel we shall require the following results.\\

In the equation (\ref{eq1.6}), substitute  $\zeta=1$ and replace $x$ with $nx$, we get
\begin{eqnarray*}
\sum_{i=0}^{\infty}p_{i}(nx)&=& A(1)e^{nx}+B(1)e^{-nx}.
\end{eqnarray*}

By taking the first four derivatives of the equation (\ref{eq1.6}) and substituting $\zeta=1$, and then replacing $x$ with $nx$, we obtained
\begin{eqnarray*}
\sum_{i=0}^{\infty}i p_{i}(nx)&=& nx(A(1)e^{nx}-B(1)e^{-nx})+A^{\prime}(1)e^{nx}+B^{\prime}(1)e^{-nx},\\
\sum_{i=0}^{\infty}i(i-1)p_{i}(nx)&=& x^2n^2(A(1)e^{nx}+B(1)e^{-nx})+2nx(A^{\prime}(1)e^{nx}-B^{\prime}(1)e^{-nx})\\
&&+(A^{\prime \prime}(1)e^{nx}+B^{\prime \prime}(1)e^{-nx})
\end{eqnarray*}
\noindent $\displaystyle
\sum_{i=0}^{\infty}i(i-1)(i-2)p_{i}(nx)$
\vskip-1mm
\begin{eqnarray*}
&=& n^3x^3(A(1)e^{nx}-B(1)e^{-nx})+3x^2n^2(A^{\prime}(1)e^{nx}+B^{\prime}(1)e^{-nx})\\
&&+3nx(A^{\prime \prime}(1)e^{nx}-B^{\prime \prime}(1)e^{-nx})+(A^{\prime \prime \prime}(1)e^{-nx}+B^{\prime \prime \prime}(1)e^{-nx}),\\
\end{eqnarray*}

\noindent $\displaystyle\sum_{i=0}^{\infty}i(i-1)(i-2)(i-3)p_{i}(nx)$
\vskip-1mm
\begin{eqnarray*}
&=& n^4x^4(A(1)e^{nx}+B(1)e^{-nx})+4n^3x^3(A^{\prime}(1)e^{nx}+B^{\prime}(1)e^{-nx})\\
&&+6n^2x^2(A^{\prime \prime}(1)e^{nx}+B^{\prime \prime}(1)e^{-nx})+4nx(A^{\prime \prime \prime}(1)e^{nx}-B^{\prime \prime \prime}(1)e^{-nx})\\
&&+(A^{\prime \prime \prime \prime}(1)e^{nx}+B^{\prime \prime \prime \prime}(1)e^{-nx}).
\end{eqnarray*}

\begin{lem}\label{lemma2.1} \cite{SVSS} The operators $T_{n}(\zeta^m;x),\,\, m=\overline{0,4}$, we have
\begin{eqnarray*}
T_{n}(1;x)&=&1;\\
T_{n}(\zeta;x)&=&\left(\frac{A(1)e^{nx}-B(1)e^{-nx}}{A(1)e^{nx}+B(1)e^{-nx}}\right)x+\frac{A^{\prime}(1)e^{nx}+B^{\prime}(1)e^{-nx}}{n[A(1)e^{nx}+B(1)e^{-nx}]};\\
T_{n}(\zeta^2;x)&=&x^2+\left(\frac{[2A^{\prime}(1)+A(1)]e^{nx}-[2B^{\prime}(1)+B(1)]e^{-nx}}{n[A(1)e^{nx}+B(1)e^{-nx}]}\right)x\\
&&+\frac{[A^{\prime \prime}(1)+A^{\prime}(1)]e^{nx}+[B^{\prime \prime}(1)+B^{\prime}(1)]e^{-nx}}{n^2[A(1)e^{nx}+B(1)e^{-nx}]};\\
T_{n}(\zeta^3;x)&=&\left(\frac{A(1)e^{nx}-B(1)e^{-nx}}{A(1)e^{nx}+B(1)e^{-nx}}\right)x^3+\frac{1}{n[A(1)e^{nx}+B(1)e^{-nx}]}\left[(3A(1)\right.\\
&&\left.+2A^{\prime}(1))e^{nx}+(3B(1)+2B^{\prime}(1))e^{-nx}\right]x^2\\
&&+\frac{1}{n^2[A(1)e^{nx}+B(1)e^{-nx}]}\left\{(A(1)+(n+6)A^{\prime}(1)+2nA^{\prime \prime})e^{nx}\right.\\
&&\left.-(B(1)+(6-n)B^{\prime}(1)+2nB^{\prime \prime}(1))e^{-nx}\right\}x\\
&&+\frac{(A^{\prime}(1)+A^{\prime \prime \prime})e^{nx}+(B^{\prime}(1)+B^{\prime \prime \prime}(1))e^{-nx}}{n^3[A(1)e^{nx}+B(1)e^{-nx}]};\\
T_{n}(\zeta^4;x)&=&x^4+\frac{\{3(2A(1)+A^{\prime}(1))e^{nx}-3(2B(1)+B^{\prime}(1))e^{-nx}\}}{n[A(1)e^{nx}+B(1)e^{-nx}]}x^3\\
&&+\frac{1}{n^2[A(1)e^{nx}+B(1)e^{-nx}]}\left\{(7A(1)+13A^{\prime}(1)+4A^{\prime \prime}(1))e^{nx}\right.\\
&&\left.+(7B(1)+11B^{\prime}(1)+4B^{\prime \prime}(1))e^{-nx}\right\}x^2\\
&&+\frac{1}{n^4[A(1)e^{nx}+B(1)e^{-nx}]}\left\{((-6-5n)A(1)+n(2(7+3n)A^{\prime}(1)\right.\\
&&\left.+2(1+6n)A^{\prime \prime}(1)+(2+n)A^{\prime \prime \prime}(1)))e^{nx}+((-6+5n)B(1)\right.\\
&&\left.+n(2(-7+3n)B^{\prime}(1)+(2-12n)B^{\prime \prime}(1)-3B^{\prime \prime \prime}(1)))e^{-nx}\right\}x\\
&&+\frac{1}{n^5[A(1)e^{nx}+B(1)e^{-nx}]}\left\{((6-5n)A^{\prime}(1)+n(13A^{\prime \prime}(1)+7A^{\prime \prime \prime}(1)\right.\\
&&\left.+A^{\prime \prime \prime \prime}(1)))e^{nx}+((6-5n)B^{\prime}(1)+n(B^{\prime \prime}(1)+5B^{\prime \prime \prime}(1)\right.\\
&&\left.+B^{\prime \prime \prime \prime}(1)))e^{-nx}\right\}.
\end{eqnarray*}
\end{lem}

\begin{lem}\label{lemma2.2} For the proposed operators $S_{n}^{\rho}(\zeta^m;x),\,\, m=\overline{0,4}$, we have
\begin{eqnarray*}
S_{n}^{\rho}(1;x)&=&1;\\
S_{n}^{\rho}(\zeta;x)&=&\frac{1}{(n\rho-c)[A(1)e^{nx}+B(1)e^{-nx}]}\left\{n\rho x[A(1)e^{nx}-B(1)e^{-nx}]\right.\\
&&\left.+\rho[A^{\prime}(1)e^{nx}+B^{\prime}(1)e^{-nx}]\right\};\\
S_{n}^{\rho}(\zeta^{2};x)&=&\frac{n^2\rho^2}{(n\rho-c)(n\rho-2c)}x^2\\
&&+\frac{1}{(n\rho-c)(n\rho-2c)[A(1)e^{nx}+B(1)e^{-nx}]}\left\{(n\rho((A(1)\right.\\
&&\left.+2\rho A^{\prime}(1))e^{nx}-(B(1)+2\rho B^{\prime}(1))e^{-nx}))x+(\rho((A^{\prime}(1)\right.\\
&&\left.+\rho A^{\prime \prime}(1))e^{nx}+(B^{\prime}(1)+\rho B^{\prime \prime}(1))e^{-nx}))\right\};\\
S_{n}^{\rho}(\zeta^{3};x)&=&\frac{1}{(n\rho-c)(n\rho-2c)(n\rho-3c)[A(1)e^{nx}+B(1)e^{-nx}]}\left\{n^3 \rho^3 (A(1)e^{nx}\right.\\
&&\left.-B(1)e^{-nx})x^3+3n^2 \rho^2 ((A(1)+\rho A^{\prime}(1))e^{nx}+(B(1)\right.\\
&&\left.+\rho B^{\prime}(1))e^{-nx})x^2+(n\rho((2A(1)+3\rho(2A^{\prime}(1)+\rho A^{\prime \prime}(1)))e^{nx}\right.\\
&&\left.-(2B(1)+3\rho(2B^{\prime}(1)+\rho B^{\prime \prime})(1))e^{-nx}))x
+\rho((2A^{\prime}(1)+\rho(3A^{\prime \prime}(1)\right.\\
&&\left.+\rho A^{\prime \prime \prime}(1)))e^{nx}+(2B^{\prime}(1)+\rho(3B^{\prime \prime}(1)+\rho B^{\prime \prime \prime}(1)))e^{-nx})\right\};\\
S_{n}^{\rho}(\zeta^{4};x)&=&\frac{n^4 \rho^4}{(n\rho-c)(n\rho-2c)(n\rho-3c)(n\rho-4c)}x^4\\
&&+\frac{1}{(n\rho-c)(n\rho-2c)(n\rho-3c)(n\rho-4c)[A(1)e^{nx}+B(1)e^{-nx}]}\times \\
&&\left\{2n^3\rho^3((3A(1)+2\rho A^{\prime}(1))e^{nx}-(3B(1)+2\rho B^{\prime}(1))e^{-nx})x^3\right.\\
&&\left.+n^2\rho^2((11A(1)+6\rho (3A^{\prime}(1)+\rho A^{\prime \prime}(1)))e^{nx}+(11B(1))\right.\\
&&\left.+6\rho(3B^{\prime}(1)+\rho B^{\prime \prime}(1)))e^{-nx})x^2+(2n\rho((3A(1)+\rho(11A^{\prime}(1\right.\\
&&\left.+9\rho A^{\prime \prime}(1)+2\rho^2 A^{\prime \prime \prime}(1)))e^{nx})-(3B(1)+\rho(11B^{\prime}(1)+9\rho B^{\prime \prime}(1)\right.\\
&&\left.+2\rho^2 B^{\prime \prime \prime}(1)))e^{-nx })x+\rho((6A^{\prime}(1)+\rho(11A^{\prime \prime}(1)+6\rho A^{\prime \prime \prime}(1)\right.\\
&&\left.+\rho^2 A^{\prime \prime \prime \prime}(1)))e^{nx}+(6B^{\prime}(1)+\rho(11B^{\prime \prime}(1)\right.\\
&&\left.+6\rho B^{\prime \prime \prime}(1)+\rho^2 B^{\prime \prime \prime \prime}(1)))e^{-nx})\right\}
\end{eqnarray*}
\end{lem}
\begin{proof}
The proof follows the lemma \ref{lemma2.1} and the equation (\ref{eq1.8}).
\end{proof}
\begin{lem}\label{lemma2.3} The central moments of the operators $S_{n}^{\rho}((\zeta-x)^m;x),$ for $m=1,\,2$ and $4$, we have
\begin{eqnarray*}
S_{n}^{\rho}(\zeta-x;x)&=&\frac{\{-2n\rho B(1)e^{-nx}+c[A(1)e^{nx}+B(1)e^{-nx}]\}}{(n\rho-c)[A(1)e^{nx}+B(1)e^{-nx}]}x\\
&&+\frac{\rho[A^{\prime}(1)e^{nx}+B^{\prime}(1)e^{-nx}]}{(n\rho-c)[A(1)e^{nx}+B(1)e^{-nx}]};\\
S_{n}^{\rho}((\zeta-x)^2;x)&=&\frac{1}{(n\rho-c)(n\rho-2c)[A(1)e^{nx}+B(1)e^{-nx}]}\left\{cn\rho((A(1)\right.\\
&&\left.+2c^2A(1))e^{nx}+(-7B(1)+(4n^2\rho^2+2c^2)B(1))e^{-nx})x^2\right.\\
&&\left.+\rho((nA(1)+4cA^{\prime}(1))e^{nx}+(4cB^{\prime}(1)-n(B(1)\right.\\
&&\left.+4\rho B^{\prime}(1)))e^{-nx})x+\rho((A^{\prime}(1)+\rho A^{\prime \prime}(1))e^{nx}+(B^{\prime}(1)\right.\\
&&\left.+\rho B^{\prime \prime}(1))e^{-nx})\right\};\\
S_{n}^{\rho}((\zeta-x)^4;x)&=& \frac{1}{(n\rho-c)(n\rho-2c)(n\rho-3c)(n\rho-4c)[A(1)e^{nx}+B(1)e^{-nx}]}\\
&&\left\{((46c^3 n\rho A(1)+24c^4 A(1)+3c^2n^2\rho^2 A(1))e^{nx}+(-146c^3 n \rho B(1)\right.\\
&&\left.-104cn^3\rho^3 B(1)+16 n^4 \rho^4 B(1)+24 c^4 B(1)\right.\\
&&\left.+211c^2 n^2 \rho^2 B(1))e^{-nx})x^4+((96\rho c^3 A^{\prime}(1)+4c^2n(9A(1)\right.\\
&&\left.+5\rho A^{\prime}(1))+3cn^2\rho A(1))e^{nx}+(96c^3 B^{\prime}(1)-36c^2nB(1)\right.\\
&&\left.-124c^2 n\rho B^{\prime}(1)-4n^3 \rho^2(3B(1)+4\rho B^{\prime}(1))+15cn^2 \rho B(1)\right.\\
&&\left.+84cn^2 \rho^2 B^{\prime}(1))e^{-nx})x^3+(\rho((72c^2(A^{\prime}(1)+\rho A^{\prime \prime}(1))\right.\\
&&\left.+3n^2 \rho A(1)-2cn(16 A(1)+3\rho(9A^{\prime}(1)+\rho A^{\prime \prime}(1)))))e^{nx}\right.\\
&&\left.+(72c^2B^{\prime}(1)+72c^2\rho B^{\prime \prime}(1)+n^2\rho (19B(1)+24\rho (2B^{\prime}(1)\right.\\
&&\left.+\rho B^{\prime \prime}(1)))+2cn(-16B(1))-3\rho(23B^{\prime}(1)+15\rho B^{\prime \prime}(1)))e^{-nx})x^2\right.\\
&&\left.-((2\rho(-8c(2B^{\prime}(1)+3\rho B^{\prime \prime}(1))-8cn^3\rho^2 B^{\prime \prime \prime}(1)+2n^4 \rho^{3}B^{\prime \prime \prime}(1)\right.\\
&&\left.+n(3B(1)+\rho(15B^{\prime}(1)+\rho(15B^{\prime \prime}(1)+2\rho B^{\prime \prime \prime}(1))))))e^{-nx}\right.\\
&&\left.-(3n A(1)+16cA^{\prime}(1)+n\rho (7A^{\prime}(1)+3\rho A^{\prime \prime}(1)))e^{nx})x\right.\\
&&\left.+((\rho(6A^{\prime}(1)+\rho(11 A^{\prime \prime}(1)+6\rho A^{\prime \prime \prime}(1)+\rho^2 A^{\prime \prime \prime \prime}(1)))e^{nx}\right.\\
&&\left.+(6 B^{\prime}(1)+\rho(11B^{\prime \prime}(1)\right.\\
&&\left.+6n^3 \rho B^{\prime \prime \prime}(1)+\rho^2 B^{\prime \prime \prime\prime}(1)))e^{-nx}))\right\}.
\end{eqnarray*}
\end{lem}
\begin{proof} The Proof of the above lemma can be easily found by using the following equalities:
\begin{eqnarray*}
S_{n}^{\rho}((\zeta-x);x)&=& S_{n}^{\rho}(\zeta;x)-xS_{n}^{\rho}(1;x);\\
S_{n}^{\rho}((\zeta-x)^2;x)&=& S_{n}^{\rho}(\zeta^2;x)-2xS_{n}^{\rho}(\zeta;x)+x^2S_{n}^{\rho}(1;x);\\
S_{n}^{\rho}((\zeta-x)^4;x)&=& S_{n}^{\rho}(\zeta^4;x)-4xS_{n}^{\rho}(\zeta^3;x)+6x^2S_{n}^{\rho}(\zeta^2;x)-4x^3S_{n}^{\rho}(\zeta;x)\\
&&+x^4S_{n}^{\rho}(1;x).
\end{eqnarray*}
\end{proof}

\begin{lem}\label{lemma2.4} The limiting values of the Lemma \ref{lemma2.3}, we have
\begin{eqnarray*}
\lim_{n\rightarrow \infty}n S_{n}^{\rho}((\zeta-x);x)&=& \frac{cx}{\rho}+\frac{A^{\prime}(1)}{A(1)};\\
\lim_{n\rightarrow \infty}n S_{n}^{\rho}((\zeta-x)^2;x)&=&\frac{x(1+cx)}{\rho};\\
\lim_{n^4\rightarrow \infty}n^2 S_{n}^{\rho}((\zeta-x)^4;x)&=& \frac{3x^2(1+cx)^2}{\rho^2}.
\end{eqnarray*}
\end{lem}

\section{Approximation Results}
In the whole paper, we take it $\varrho_{n}(x)=S_{n}^{\rho}((\zeta-x)^2;x)$.
In order to provide continuous functions for positive linear operators in a compact space which strongly converges to the identity operators, Bohman-Korovkin delegated a simple and convincing theorem. Using the Bohman-Korovkin theorem, the following theorem provides an approximation result for the operators $S_{n}^{\rho}$.

\begin{thm}\label{theorem3.1}
Let $f$ be a continuous function on the interval $[0,\infty)$, then
\begin{eqnarray*}
\lim_{n\rightarrow \infty}S_{n}^{\rho}(f;x)=f(x)
\end{eqnarray*}holds uniformly in $x\in [0,\infty)$.
\end{thm}

\begin{proof}
The proof is based on the conditions by the Korovkin's theorem. By using Lemma \ref{lemma2.2}, we obtain
$$\lim_{n\rightarrow \infty}S_{n}^{\rho}(\zeta^{r};x)=x^{r},\,\,r={0,1,2}.$$
\end{proof}
The rate of convergence of the given operators $S^{\rho}_{n}$  for the Lipschitz space is shown in the results below. We get a pointwise approximation because $x$ is contained in the denominator. $x$ drop first to the uniform convergence from the Sz$\acute{a}$sz operators \cite{OZ}.
The Lipschitz class $f\in Lip_{M}^{*}(\alpha)$ holds, if the inequality
\begin{eqnarray*}
|f(\zeta)-f(x)|\leq M|\zeta-x|^{\alpha},\,\,\text{for each}\,\,\, \zeta,x\in(0,1),
\end{eqnarray*}holds. For $1<p<\infty$ and $\frac{1}{p}+\frac{1}{q}=1$, the H\"{o}lder inequality is
$$\sum_{i=0}^{\tau}|\xi_{i}\eta_{i}|\leq \left(\sum_{i=0}^{\tau}(\xi_{i})^{p}\right)^{\frac{1}{p}}\left(\sum_{i=0}^{\tau}(\eta_{i})^{q}\right)^{\frac{1}{q}},$$
where $\xi_{i},\eta_{i}\in \mathbb{K}$ (real field $\mathbb{K}$).

\begin{thm}\label{theorem3.2}
Let $f\in Lip^{*}_{M}(\alpha)$ and $\alpha \in (0,1]$. Then the inequality holds:
\begin{eqnarray*}
|S_{n}^{\rho}(f;x)-f(x)|\leq M\left(\frac{\varrho_{n}(x)}{x}\right)^{\frac{\alpha}{2}},
\end{eqnarray*}for all $x\in (0,\infty)$.
\end{thm}
\begin{proof}
In light of the proposed operators linearity and positivity (\ref{eq1.8}), we write as
\begin{eqnarray*}
|S_{n}^{\rho}(f;x)-f(x)|\leq \int_{0}^{\infty}A_{n}(x,\zeta)|f(\zeta)-f(x)|d\zeta.
\end{eqnarray*}Applying H$\ddot{o}$lder's inequality with $p=\frac{2}{\alpha}$ and $q=\frac{2}{2-\alpha}$ and Lemma \ref{lemma2.2}, we obtain
\begin{eqnarray*}
|S_{n}^{\rho}(f;x)-f(x)|&\leq& \left(\int_{0}^{\infty}A_{n}(x,\zeta)|f(\zeta)-f(x)|^{\frac{2}{\alpha}}d\zeta\right)^{\frac{\alpha}{2}}\left(\int_{0}^{\infty}A_{n}(x,\zeta)d\zeta\right)^{\frac{2-\alpha}{2}}\\
&\leq& \left(\int_{0}^{\infty}A_{n}(x,\zeta)|f(\zeta)-f(x)|^{\frac{2}{\alpha}}d\zeta\right)^{\frac{\alpha}{2}}\\
&\leq& M\left(\int_{0}^{\infty}A_{n}(x,\zeta)\frac{(\zeta-x)^2}{(\zeta+x)}d\zeta\right)^{\frac{\alpha}{2}}\\
&\leq& M\left(\frac{\varrho_{n}(x)}{x}\right)^{\frac{\alpha}{2}}.
\end{eqnarray*}We get the required result.
\end{proof}
In the present theorem, we set up a Voronvaskaja type approximation theorem.
\begin{thm}\label{theorem3.3}
Let $f\in C_{2}[0,\infty)$ and $f^{\prime},f^{\prime \prime}$ exists  at a pont $x\in[0,\infty)$, then the equality holds
\begin{eqnarray*}
\lim_{n\rightarrow \infty}n(S_{n}^{\rho}(f;x)-f(x))&=&\left(\frac{cx}{\rho}+\frac{A^{\prime}(1)}{A(1)}\right)f^{\prime}(x)\\
&&+\left(\frac{cx^2}{\rho}+\frac{x(A(1)+(2-A^{\prime}(1)\rho))}{A(1)\rho}\right)\frac{f^{\prime \prime}(x)}{2}.
\end{eqnarray*}If $f^{\prime \prime}$ is continuous on $[0,\infty)$, then the above result holds uniformly in the interval $[0,b]\subset [0,\infty)$ with $b>0.$
\end{thm}

\begin{proof}
We know that the Taylor's series expansion
\begin{eqnarray}\label{eq3.1}
f(\zeta)=f(x)+f^{\prime}(x)(\zeta-x)+\frac{f^{\prime \prime}(x)}{2}(\zeta-x)^2+\sigma(\zeta,x)(\zeta-x)^2,
\end{eqnarray}where $\sigma(\zeta,x)\in C_{2}[0,\infty)$ and $\displaystyle\lim_{n\rightarrow \infty}\sigma(\zeta,x)=0$. Applying the given operator $S_{n}^{\rho}$ from the equation (\ref{eq1.8}) on both sides of the equation (\ref{eq3.1}), we  have
\begin{eqnarray}\label{eq3.2}
\lim_{n\rightarrow \infty}n(S_{n}^{\rho}(f;x)-f(x))&=&\lim_{n\rightarrow \infty}n S_{n}^{\rho}((\zeta-x);x)f^{\prime}(x)+\lim_{n\rightarrow \infty}n S_{n}^{\rho}((\zeta-x)^2;x)\frac{f^{\prime \prime}(x)}{2}\nonumber\\
&&+\lim_{n\rightarrow \infty}n S_{n}^{\rho}(\sigma(\zeta,x)(\zeta-x)^2;x).
\end{eqnarray}
In the last term of the equation (\ref{eq3.2}), using the Cauchy-Schwarz inequality, we obtain
$$n S_{n}^{\rho}(\sigma(\zeta,x)(\zeta-x)^2;x)\leq \sqrt{S_{n}^{\rho}(\sigma^{2}(\zeta,x);x)}\sqrt{n^2S_{n}^{\rho}((\zeta-x)^4;x)}.$$
Since $\sigma(\zeta,x)\rightarrow 0,$ as $t\rightarrow x$, applying Theorem \ref{theorem3.1}, we get $\displaystyle \lim_{n\rightarrow \infty} S_{n}^{\rho}(\sigma^2(\zeta,x);x)=\sigma^2(x;x)=0$. And applying Lemma \ref{lemma2.4}, for large $n$, and $x\in [0,\infty)$, we have
\begin{eqnarray}\label{eq3.3}
n^2S_{n}^{\rho}((\zeta-x)^4;x)=O(1).
\end{eqnarray}Hence,
\begin{eqnarray}\label{eq3.4}
\lim_{n\rightarrow \infty}nS_{n}^{\rho}(\sigma(\zeta,x)(\zeta-x)^2;x)=0.
\end{eqnarray}The required result are obtained from the equations (\ref{eq3.2}), (\ref{eq3.4}) and Lemma \ref{lemma2.4}, we have
\begin{eqnarray*}
\lim_{n\rightarrow \infty}n(S_{n}^{\rho}(f;x)-f(x))&=&\left(\frac{cx}{\rho}+\frac{A^{\prime}(1)}{A(1)}\right)f^{\prime}(x)\\
&&+\left(\frac{cx^2}{\rho}+\frac{x(A(1)+(2-A^{\prime}(1)\rho))}{A(1)\rho}\right)\frac{f^{\prime \prime}(x)}{2}.
\end{eqnarray*}
\end{proof}

The uniformity declaration is necessarily dependent on the regular continuity of $f^{\prime \prime}$ on $[0,b]$, and other conclusions are true in general if $[0, b],$ where $b>0$. Using the classical modulus of continuity, we analysed the proposed operators approximation results in the subsequent theorem.
%
\begin{thm}\label{theorem3.4}
For $f\in C_{2}[0,\infty)$, then the inequality holds
\begin{eqnarray*}
|S_{n}^{\rho}(f;x)-f(x)|\leq 4 M_{f}(1+x^2)\varrho_{n}(x)+2\omega_{b+1}(f;\sqrt{\varrho_{n}(x)}).
\end{eqnarray*}where $\omega(f;\varrho_{n}(x))$ is the modulus of continuity of $f$ on $[0,b+1]$.
\end{thm}
\begin{proof}
From \cite{EIEAG}, for $\zeta\in (b+1,\infty)$ and $x\in [0,b]$, we have
\begin{eqnarray*}
|f(\zeta)-f(x)|\leq 4M_{f}(t-x)^2(1+x^2)+\left(1+\frac{|\zeta-x|}{\delta}\right)\omega_{b+1}(f,\delta),\,\,\delta>0.
\end{eqnarray*}Applying the cauchy-Schwarz inequality, we get
\begin{eqnarray*}
|S_{n}^{\rho}(f;x)-f(x)|&\leq& 4M_{f}(1+x^2)S_{n}^{\rho}((\zeta-x)^2;x)\\
&&+\left(1+\frac{(S_{n}^{\rho}((\zeta-x)^2;x))^{\frac{1}{2}}}{\delta}\right)\omega_{b+1}(f,\delta)\\
&\leq& 4 M_{f}(1+x^2)\varrho_{n}(x)+\omega_{b+1}(f,\delta)\left(1+\frac{\sqrt{\varrho_{n}(x)}}{\delta}\right).
\end{eqnarray*} Hence suppose $\delta=\sqrt{\varrho_{n}(x)}$, we obtained the result.
\end{proof}

In \cite{DZTV}, Ditzian-Totik and Peetre's K-functional initiated the modulus of continuity. They considered the function $\psi^2(x)=x(1+x)$ and $f\in C_{B}[0,\infty)$ the space of all bounded and continuous functions on $[0,\infty)$ endowed with the norm $\displaystyle\|f\|=\sup_{x\in[0,\infty)}|f(x)|$. The modulus $\omega^{*}_{\psi^{\tau}}$, where $0\leq \tau\leq 1$ is defined as
\begin{eqnarray*}
\omega^{*}_{\psi^{\tau}}(f,\zeta)=\sup_{0\leq h\leq t}\sup_{x\pm\frac{h\psi^{\tau}(x)}{2}\in[0,\infty)}\left|f\left(x+\frac{h\psi^{\tau}(x)}{2}\right)-f\left(x-\frac{h\psi^{\tau}(x)}{2}\right)\right|,
\end{eqnarray*}and the appropriate $k-$functional is given by
$$K^{*}_{\psi^{\tau}}(f,\zeta)=\inf_{g\in W_{\tau}}\{\|f-g\|+\zeta\|\psi^{\tau}g^{\prime}\|\},$$
where $W_{\tau}=\{g:g\in AC_{loc}[0,\infty):\|\psi^{\tau}g^{\prime}\|<\infty\}$, $AC_{loc}$ denotes the space of locally absolutely continuous function on $[0,\infty)$. For a constant $M>0$ such that
\begin{eqnarray}\label{eq3.5}
M^{-1}\omega^{*}_{\psi^{\tau}}(f,\zeta)\leq K^{*}_{\psi^{\tau}}(f,\zeta)\leq M \omega^{*}_{\psi^{\tau}}(f,\zeta).
\end{eqnarray}
\begin{thm}\label{theorem3.5}
For $f\in C_{B}[0,\infty)$, then the inequality holds for the large value of $n$
\begin{eqnarray*}
|S_{n}^{\rho}(f;x)-f(x)|\leq C\omega^{*}_{\psi^{\tau}}\left(f;\frac{\psi^{1-\tau}(x)}{\sqrt{n}}\right),
\end{eqnarray*}where $C$ is a constant independent of $f$ and $n$.
\end{thm}

\begin{proof}
By the definition of $K-$functional, choose $g\in W_{\tau}$ such that
\begin{eqnarray}\label{eq3.6}
\|f-g\|+\frac{\psi^{1-\tau}(x)}{\sqrt{n}}\|\psi^{\tau}g^{\prime}\|\leq 2 K^{*}_{\psi^{\tau}}\left(f;\frac{\psi^{1-\tau}(x)}{\sqrt{n}}\right).
\end{eqnarray} Now,
\begin{eqnarray}\label{eq3.7}
|S^{\rho}_{n}(f;x)-f(x)|&\leq& |S^{\rho}_{n}(f-g;x)|+|S^{\rho}_{n}(g;x)-g(x)|+|g(x)-f(x)|\nonumber\\
&\leq& 2\|f-g\|+|S^{\rho}_{n}(g;x)-g(x)|.
\end{eqnarray}Consider
$$g(t)=g(x)+\int_{x}^{\zeta}g^{\prime}(u)du$$
and
\begin{eqnarray}\label{eq3.8}
|S^{\rho}_{n}(g;x)-g(x)|\leq S_{n}^{\rho}\left(\left|\int_{x}^{\zeta}g^{\prime}(u)du\right|;x\right).
\end{eqnarray}Applying H$\ddot{o}$lder's inequality, we have
\begin{eqnarray*}
\left|\int_{x}^{t}g^{\prime}(u)du\right|\leq \|\psi^{\tau}g^{\prime}\|\left|\int_{x}^{\zeta}\frac{du}{\psi^{\tau}(u)}\right|\leq \|\psi^{\tau}g^{\prime}\||\zeta-x|^{1-\tau}\left|\int_{x}^{\zeta}\frac{du}{\psi(u)}\right|^{\tau},
\end{eqnarray*}we may write
\begin{eqnarray*}
\left|\int_{x}^{\zeta}\frac{du}{\psi(u)}\right|\leq \left|\int_{x}^{\zeta}\frac{du}{\sqrt{u}}\right|\left(\frac{1}{\sqrt{1+x}}+\frac{1}{\sqrt{1+\zeta}}\right).
\end{eqnarray*}
Using the inequality $|\alpha+\beta|^{r}\leq |\alpha|^r+|\beta|^r,\,\,0\leq r\leq 1$, we have
\begin{eqnarray}\label{eq3.9}
\left|\int_{x}^{\zeta}g^{\prime}(u)du\right|&\leq& \frac{2^{\tau}\|\psi^{\tau}g^{\prime}\||\zeta-x|}{x^{\frac{\tau}{2}}}\left(\frac{1}{\sqrt{1+x}}+\frac{1}{\sqrt{1+\zeta}}\right)^{\tau}\nonumber\\
&\leq& \frac{2^{\tau}\|\psi^{\tau}g^{\prime}\||\zeta-x|}{x^{\frac{\tau}{2}}}\left(\frac{1}{(1+x)^{\frac{\tau}{2}}}+\frac{1}{(1+\zeta)^{\frac{\tau}{2}}}\right).
\end{eqnarray} Thus, from the equation (\ref{eq3.8}), (\ref{eq3.9}), using Cauchy-Schwartz inequality, Theorem \ref{theorem3.1}, and sufficiently large $n$, we obtain
\begin{eqnarray}\label{eq3.10}
|S_{n}^{\rho}(g;x)-g(x)|&\leq& \frac{2^{\tau}\|\psi^{\tau}g^{\prime}\|}{x^{\frac{\tau}{2}}}S_{n}^{\rho}\left(|\zeta-x|\left(\frac{1}{(1+x)^{\frac{\tau}{2}}}+\frac{1}{(1+\zeta)^{\frac{\tau}{2}}}\right);x\right)\nonumber\\
&\leq& \frac{2^{\tau}\|\psi^{\tau}g^{\prime}\|}{x^{\frac{\tau}{2}}}\left(\frac{1}{(1+x)^{\frac{\tau}{2}}}\sqrt{S_{n}^{\rho}((\zeta-x)^2;x)}\nonumber\right.\\
&&\left.+\sqrt{S_{n}^{\rho}((\zeta-x);x)}\sqrt{S_{n}^{\rho}((1+\zeta)^{-\tau};x)}\right)\nonumber\\
&\leq& 2^{\tau}\|\psi^{\tau}g^{\prime}\|\sqrt{S_{n}^{\rho}((\zeta-x)^2;x)}\left\{\psi^{-\tau}(x)+x^{\frac{-\tau}{2}}\sqrt{S_{n}^{\rho}((1+\zeta)^{-\tau};x)}\right\}\nonumber\\
&\leq& 2^{\tau}C\|\phi^{\tau}g^{\prime}\|\frac{\phi(x)}{\sqrt{n}}\{\psi^{-\tau}(x)+x^{\frac{-\tau}{2}}(1+x)^{\frac{-\tau}{2}}\}\nonumber\\
&\leq& 2^{\tau+1}\frac{\|\psi^{\tau}g^{\prime}\|\psi^{1-\tau}(x)}{\sqrt{n}}.
\end{eqnarray} Hence, from the equations (\ref{eq3.6}- \ref{eq3.8}), and (\ref{eq3.10}), we get the required result
\begin{eqnarray*}
|S_{n}^{\rho}(f;x)-f(x)|&\leq& 2\|f-g\|+2^{\tau+1}C\|\psi^{\tau}g^{\prime}\|\frac{\psi^{1-\tau}}{\sqrt{n}}\\
&\leq& C\left\{\|f-g\|+\frac{\psi^{1-\tau}}{\sqrt{n}}\|\psi^{\tau}g^{\prime}\|\right\}\leq 2CK^{*}_{\psi^{\tau}}\left(f;\frac{\psi^{1-\tau}(x)}{\sqrt{n}}\right)\\
&\leq& C \omega^{*}_{\psi^{\tau}}\left(f;\frac{\psi^{1-\tau}(x)}{\sqrt{n}}\right).
\end{eqnarray*}
\end{proof}
Here, we'll discuss about the consequences of the Ras'as and Steklov functions and determine approximation results by using the second-order continuity modulus. The modulus of continuity is defined as: $$\omega_{2}(\phi,\delta)= \sup_{0<\zeta \leq \delta}\|\phi(.+2\zeta)-2\phi(.+\zeta)+\phi(.)\|.$$
Now consider $\{L_{n}\}_{n\geq 0}$ be a  sequence of linear positive operators with virtue $L_{n}(e_{i};x)=x^{i}$. Then according to the Ras\'{a}s result \cite{GIRI}, we have
\begin{eqnarray*}
|L_{n}(g;x)-g(x)|&\leq& \|g^{\prime}\|\sqrt{L_{n}((\zeta-x)^2;x)}+\frac{1}{2}\|g^{\prime \prime}\|L_{n}((\zeta-x)^2;x),
\end{eqnarray*}where $g\in C^{2}[0,a]$.
And for $f\in C[a,b]$, the second-order Steklov function is as follows
\begin{eqnarray*}
f_{h}(x)=\frac{1}{h}\int_{-h}^{h}\left(1-\frac{|\zeta|}{h}\right)f(h;x+\zeta)d\zeta,\,\,x\in [a,b],
\end{eqnarray*}where $f(h;.):[a-h,b+h]\rightarrow \mathbb{R},\,h>0$ by
$$
f_{h}(x)=
\begin{cases}
P_{-}(x);\,\,\,\,a-h\leq x\leq a\\
f(x);\,\,\,\,\,\,\,\,\,a\leq x\leq b\\
P_{+}(x);\,\,\,\,b<x\leq b+h
\end{cases}
$$
and $P_{-},\,P_{+}$ they are linear best approximation to the function $f$ on the given interval.

\begin{thm}\label{theorem3.5}
Let $\phi\in C[0,\infty)$. Then, we obtain
\begin{eqnarray*}
 |S_{n}^{\rho}(\phi;x)-\phi(x)|&\leq& \left(\frac{3}{2}+\frac{3a}{4}+\frac{3h^2}{4}\right)\omega_{2}(\phi;h)+\frac{2h^2}{a}\|\phi\|.
\end{eqnarray*}
\end{thm}

\begin{proof}
From some calculations and using well-known properties, we have
\begin{eqnarray}\label{eq3.11}
|S_{n}^{\rho}(\phi;x)-\phi(x)|&\leq& S_{n}^{\rho}(|\phi-\phi_{h}|;x)+|S_{n}^{\rho}(\phi_{h};x)-\phi_{h}(x)|+|\phi_{h}(x)-\phi(x)|\nonumber\\
&\leq& 2\|\phi-\phi_{h}\|+|S_{n}^{\rho}(\phi_{h};x)-\phi_{h}(x)|,
\end{eqnarray} where, $\phi_{h}\in C^{2}[0,a)$ be the second-order Steklov function of $\phi$. From the Rasa's result and Landau inequality, we have
\begin{eqnarray}\label{eq3.12}
|S_{n}^{\rho}(\phi;x)-\phi(x)|&\leq& \|\phi^{\prime}_{h}\|\sqrt{S_{n}^{\rho}\left((\zeta-x)^2;x\right)}+\frac{1}{2}\|\phi^{\prime \prime}_{h}\|S_{n}^{\rho}\left((\zeta-x)^2;x\right)\nonumber\\
&\leq& \left(\frac{2}{a}\|\phi_{h}\|+\frac{a}{2}\|\phi^{\prime \prime}_{h}\|\right)\sqrt{S_{n}^{\rho}\left((\zeta-x)^2;x\right)}\nonumber\\
&&+\frac{1}{2}\|\phi^{\prime \prime}_{h}\|S_{n}^{\rho}\left((\zeta-x)^2;x\right)\nonumber\\
&\leq& \left(\frac{2}{a}\|\phi\|+\frac{3a}{4h^2}\omega_{2}(\phi;h)\right)\sqrt{S_{n}^{\rho}\left((\zeta-x)^2;x\right)}\nonumber\\
&&+\frac{3}{4h^2}\omega_{2}(\phi;h)S_{n}^{\rho}\left((\zeta-x)^2;x\right).
\end{eqnarray}Using the relation given by Zhuk \cite{ZVV}, between Steklov function and $\omega_{2}(\phi;h)$ as: $\|\phi-\phi_{h}\|\leq \frac{3}{4}\omega_{2}(\phi;h)$, from the equations (\ref{eq3.11}, \ref{eq3.12}), we obtain
\begin{eqnarray*}
 |S_{n}^{\rho}(\phi;x)-\phi(x)|&\leq& \frac{3}{2}\omega_{2}(\phi;h)+\left(\frac{2}{a}\|\phi\|+\frac{3a}{4h^2}\omega_{2}(\phi;h)\right)\sqrt{S_{n}^{\rho}\left((\zeta-x)^2;x\right)}\\
 &&+\frac{3}{4h^2}\omega_{2}(\phi;h)S_{n}^{\rho}\left((\zeta-x)^2;x\right).
\end{eqnarray*}Choose $h^2=\sqrt{S_{n}^{\rho}\left((\zeta-x)^2;x\right)}$ and, by the simple calculation, we get the required result.
\end{proof}

\section{Weighted Approximation}
Consider the space $\displaystyle C_{2}^{*}[0,\infty)=\left\{f\in C_{2}[0,\infty):\lim_{x\rightarrow \infty}\frac{|f(x)|}{1+x^2}\,\text{exists\,and\,is\,finite}\right\}$. In this section, we study the approximation results on the space $C_{2}^{*}[0,\infty)$ using the weighted modulus of continuity. We know that in general the classical modulus of continuity of first order does not tends to zero on an infinite interval. So here we consider the weighted modulus of continuity defined by Y$\ddot{u}ksel$ and Ispir $\Omega^{*}(f;\delta)$, which is defined as $$\Omega^{*}(f;\delta)=\sup_{x\in[0,\infty),\,0<h<\delta}\frac{|f(x+h)-f(x)|}{1+(x+h)^2},$$ where $f\in C_{2}^{*}[0,\infty)$.\\
Before we discuss our main results firstly we study the following Lemma.
\begin{lem}\label{lemma4.1} Let $f\in C_{2}^{*}[0,\infty)$. Then the following results hold:
\begin{itemize}
  \item $\Omega^{*}(f;\delta)$ is a monotonically increasing function of $\delta;$
  \item $\displaystyle\lim_{\delta\rightarrow 0^{+}}\Omega^{*}(f;\delta)=0;$
  \item For each $m\in \mathbb{N},\Omega^{*}(f;m\delta)\leq m\Omega^{*}(f;\delta);$
  \item For each $\lambda\in(0,\infty),\,\Omega^{*}(f;\lambda\delta)\leq (1+\lambda)\Omega^{*}(f;\delta).$
\end{itemize}
\end{lem}
Now, we develop an approximation theorem using the weighted space of the continuous function $C_{2}^{*}[0,\infty)$ of the given operators $S_{n}^{\rho}$.

\begin{thm}\label{theorem4.1}
For $f\in C_{2}^{*}[0,\infty)$ and $\alpha>0$, we have
\begin{eqnarray*}
\lim_{n\rightarrow \infty}\sup_{x\in[0,\infty)}\frac{|S_{n}^{\rho}(f;x)-f(x))|}{(1+x^2)^{1+\alpha}}=0.
\end{eqnarray*}
\end{thm}

\begin{proof}
Suppose that $x_{0}$ be an arbitrary fixed point, then
\begin{eqnarray}\label{eq4.1}
&&\sup_{x\in[0,\infty)}\frac{|S_{n}^{\rho}(f;x)-f(x))|}{(1+x^2)^{1+\alpha}}\leq\sup_{x\leq x_{0}}\frac{|S_{n}^{\rho}(f;x)-f(x))|}{(1+x^2)^{1+\alpha}}+\sup_{x>x_{0}}\frac{|S_{n}^{\rho}(f;x)-f(x))|}{(1+x^2)^{1+\alpha}}\nonumber\\
&&\leq\|S_{n}^{\rho}(f;.)-f\|_{C[0,x_{0}]}+\|f\|_{2}\sup_{x>x_{0}}\frac{S_{n}^{\rho}(1+\zeta^2);x}{(1+x^2)^{1+\alpha}}+\sup_{x>x_{0}}\frac{|f(x)|}{(1+x^2)^{1+\alpha}}.
\end{eqnarray}
Since $|f(x)|\leq \|f\|_{2}(1+x^2)$, we have
\begin{eqnarray*}
\sup_{x>x_{0}}\frac{|f(x)|}{(1+x^2)^{1+\alpha}}\leq \frac{\|f\|_{2}}{(1+x_{0}^2)^{\alpha}}.
\end{eqnarray*}Choose a number $\upsilon>0$ and $x_{0}$ to be large, then
\begin{eqnarray}\label{eq4.2}
\frac{\|f\|_{2}}{(1+x_{0}^2)^{\alpha}}<\frac{\upsilon}{6}\Rightarrow \sup_{x>x_{0}}\frac{\|f\|_{2}}{(1+x_{0}^2)^{\alpha}}<\frac{\upsilon}{6}.
\end{eqnarray}Using Theorem \ref{theorem3.1}, there exist $n_{1}\in \mathbb{N}$ such that
\begin{eqnarray*}
\|f\|_{2}\frac{S_{n}^{\rho}(1+\zeta^2;x)}{(1+x^2)^{1_\alpha}}&\leq& \frac{\|f\|_{2}}{(1+x^2)^{\alpha}}\left(1+x^2+\frac{\upsilon}{3\|f\|_{2}}\right),\,\,\forall\,n>n_{1}\\
&\leq& \frac{\|f\|_{2}}{(1+x_{0}^2)^{\alpha}}+\frac{\upsilon}{3},\,\,\forall\,n>n_{1},\,x>x_{0}.
\end{eqnarray*}Hence,
\begin{eqnarray}\label{eq4.3}
\|f\|_{2}\sup_{x>x_{0}}\frac{S_{n}^{\rho}(1+\zeta^2;x)}{(1+x^2)^{1+\alpha}}&\leq& \frac{\upsilon}{2},\,\,\forall n>n_{1}.
\end{eqnarray}Applying Theorem \ref{theorem3.3}, there exist $n_{2}\in \mathbb{N}$ such that
\begin{eqnarray}\label{eq4.4}
\|S_{n}^{\rho}(f;.)-f\|_{c[0,x_{0}]}<\frac{\upsilon}{3},\,\,n>n_{2}.
\end{eqnarray}Consider $n_{0}=max(n_{1},n_{2})$. Thus combining (\ref{eq4.1}-\ref{eq4.4}), we obtained the result
\begin{eqnarray*}
\sup_{x\in[0,\infty)}\frac{|S_{n}^{\rho}(f;x)-f(x)|}{(1+x^2)^{1+\alpha}}<\upsilon,\,\,\forall\,n>n_{0}.
\end{eqnarray*}
In the next theorem, we find the order of approximation for the weighted space corresponding to the proposed operators $S_{n}^{\rho}$.
\end{proof}

\begin{thm}\label{theorem4.2}
Let $f\in C_{2}^{*}[0,\infty)$, and sufficiently large $n$, we have
\begin{eqnarray}\label{eq4.5}
|S_{n}^{\rho}(f;x)-f(x)|\leq C(x)\Omega^{*}\left(f;\frac{1}{\sqrt{n}}\right)
\end{eqnarray}
\end{thm}
\begin{proof}
For $x\in(0,\infty)$ and $\delta>0$, using definition of $\Omega^{*}\left(f;\delta\right)$ and Lemma \ref{lemma4.1}, we have
\begin{eqnarray*}
|f(\zeta)-f(x)|&\leq& (1+(x+|x-t|)^2)\Omega^{*}(f;|t-x|)\\
&\leq& 2(1+x^2)(1+(t-x)^2)\left(1+\frac{|\zeta-x|}{\delta}\right)\Omega^{*}(f;\delta).
\end{eqnarray*} Applying $S_{n}^{\rho}(.;x)$ on both sides, we get
\begin{eqnarray}\label{eq4.6}
|S_{n}^{\rho}(f;x)-f(x)|&\leq& 2(1+x^2)\Omega^{*}(f;\delta)\left\{1+S_{n}^{\rho}((\zeta-x)^2;x)\nonumber\right.\\
&&\left.+S_{n}^{\rho}\left((1+(\zeta-x)^2\frac{|\zeta-x|}{\delta};x)\right)\right\}.
\end{eqnarray}Applying the Cauchy-Schwarz inequality of the above equation, we get
\begin{eqnarray}\label{eq4.7}
S_{n}^{\rho}\left((1+(\zeta-x)^2)\frac{|\zeta-x|}{\delta};x\right)&\leq& \frac{(S_{n}^{\rho}((\zeta-x)^2;x))^{\frac{1}{2}}}{\delta}\nonumber\\
&&+\left\{\frac{1}{\delta}(S_{n}^{\rho}((\zeta-x)^4;x))^{\frac{1}{2}}\nonumber \right.\\
&&\left.\times(S_{n}^{\rho}((\zeta-x)^2;x))^{\frac{1}{2}}\right\},
\end{eqnarray}compile the equations (\ref{eq4.5}-\ref{eq4.7}) and, taking $\delta=\frac{1}{\sqrt{n}}$, we get the required result.
\end{proof}

\section{Rate of convergence using derivative of bounded variation (DBV)}

So many researchers in the area of approximation theory achieved convergence results for the sequence of linear positive operators using the derivative of bounded variation. Let $DBV[0,\infty)$ be the class of all functions in $C_{2}[0,\infty)$ having a derivative that is local of bounded variation on $[0,\infty)$. The function $f\in DBV [0,\infty)$ is defined as $$f(x)=\int_{0}^{x}g(t)+f(0),$$ where $g$ is a function of bounded variation on every finite subintervals of $[0,\infty)$.

\begin{lem}\label{lemma5.1}\cite{MSPNSA}
Let $\theta=\theta(n)\rightarrow 0$, as $n\rightarrow \infty$ and, $\displaystyle \lim_{n\rightarrow \infty}n\theta (n)=l\in \mathbb{R}$. For adequately large $n$, we have
\begin{itemize}
\item[(i)] $\xi_{n}(x,\zeta)=\int_{0}^{\zeta}A_{n}(x,\zeta)du\leq \frac{C_{1}|\varpi(x)|}{(x-\zeta)^2}$\\
\item[(ii)] $1-\xi_{n}(x,\zeta)=\int_{\zeta}^{\infty}A_{n}(x,\zeta)du\leq \frac{C_{1}|\varpi(x)|}{(\zeta-x)^2}$,
\end{itemize} where $x\in (0,\infty)$ and, $\varpi(x)=\frac{x(1+cx)}{\rho}$.
\end{lem}
\begin{proof}
From the Lemma \ref{lemma2.2}, we have
\begin{eqnarray*}
\xi_{n}(x,\zeta)&=&\int_{0}^{\zeta}A_{n}(x,\zeta)du\\
&\leq& \int_{0}^{\zeta}\left(\frac{x-u}{x-t}\right)^2 A_{n}(x,\zeta)(x,u)du\\
&\leq& \frac{1}{(x-t)^2}S_{n}^{\rho}((u-x)^2;x)\\
&\leq& \frac{C_{1}|\varpi(x)|}{(x-\zeta)^2},
\end{eqnarray*} when $n$ is large. Similarly proof for $(ii)$.
\end{proof}

\begin{thm}\label{theorem5.1}Let $f\in DBV[0,\infty)$, for every $x\in (0,\infty)$ and, sufficiently large $n$, we have
\begin{eqnarray*}
|S_{n}^{\rho}(f;x)-f(x)|&\leq& \left\{\frac{\{-2n\rho B(1)e^{-nx}+c[A(1)e^{nx}+B(1)e^{-nx}]\}}{(n\rho-c)[A(1)e^{nx}+B(1)e^{-nx}]}x\right.\\
&&\left.+\frac{\rho[A^{\prime}(1)e^{nx}+B^{\prime}(1)e^{-nx}]}{(n\rho-c)[A(1)e^{nx}+B(1)e^{-nx}]}\right\}\left|\frac{f^{\prime}(x+)+f^{\prime}(x-)}{2}\right|\\
&&+\sqrt{C_{1}|\varpi(x)|}\left|\frac{f^{\prime}(x+)-f^{\prime}(x-)}{2}\right|+\frac{C_{1}|\varpi(x)|}{x}\sum_{i=1}^{[\sqrt{n}]}\left(\bigvee_{x-\frac{x}{i}}^{x}f^{\prime}_{x}\right)\\
&&+\frac{x}{\sqrt{n}}\left(\bigvee_{x-\frac{x}{\sqrt{n}}}^{x}f^{\prime}_{x}\right)+\left(4M_{f}+\frac{M_{f}+|f(x)|}{x^2}\right)C_{1}|\varpi(x)|\\
&&+|f^{\prime}(x+)|\sqrt{C_{1}|\varpi(x)|}\\
&&+\frac{C_{1}|\varpi(x)|}{x^2}|f(2x)-f(x)-xf^{\prime}(x+)|+\frac{x}{\sqrt{n}}\left(\bigvee_{x}^{x+\frac{x}{\sqrt{n}}}f^{\prime}_{x}\right)\\
&&+\frac{C_{1}|\varpi(x)|}{x}\sum_{i=1}^{[\sqrt{n}]}\left(\bigvee_{x}^{x+\frac{x}{\sqrt{n}}}f^{\prime}_{x}\right).
\end{eqnarray*} Here $C_{1}$ be a positive constant and $\displaystyle\bigvee_{a}^{b}f$ denotes the total variation of the function $f$ on $[a,b]$ and, $f^{\prime}_{x}$ is defined as:
\begin{eqnarray}\label{eq5.1}
 f_{x}^{\prime}(\zeta)=\begin{cases}
 f^{\prime}(\zeta)-f^{\prime}(x-),\,\,\,0\leq \zeta<x,\\
 0,\,\,\,\,\,\,\,\,\,\,\,\,\,\,\,\,\,\,\,\,\,\zeta=x,\\
 f^{\prime}(\zeta)-f^{\prime}(x+),\,\,\,x<\zeta<\infty.
 \end{cases}
\end{eqnarray}
\end{thm}

\begin{proof}
Suppose for any function $f\in DBV[0,\infty)$, and using equation (\ref{eq5.1}), we write as
\begin{eqnarray}\label{eq5.2}
f^{\prime}(u)&=& \frac{1}{2}(f^{\prime}(x)+f^{\prime}(x-))+f_{x}^{\prime}(u)+\frac{1}{2}(f^{\prime}(x+)-f^{\prime}(x-))sgn(u-x)\nonumber\\
&&+\delta_{x}(u)\left(f^{\prime}(u)-\frac{1}{2}(f^{\prime}(x+)+f^{\prime}(x-))\right),
\end{eqnarray}where
\begin{eqnarray*}
 \delta_{x}(u)=\begin{cases}
 1,\,\,u=x,\\
 0,\,\,u\neq x.
  \end{cases}
\end{eqnarray*}
We have $S_{n}^{\rho}(1;x)=1$ and, using (\ref{eq5.2}) for every $x\in (0,\infty)$, we obtain
\begin{eqnarray}\label{eq5.3}
S_{n}^{\rho}(f;x)-f(x)&=& \int_{0}^{\infty}A_{n}(x,\zeta)(f(\zeta)-f(x))d{\zeta}\nonumber\\
&=& \int_{0}^{\infty}A_{n}(x,\zeta)\left(\int_{x}^{\zeta}f^{\prime}(u)\right)d(\zeta)\nonumber\\
&=&-\int_{0}^{x}\left(\int_{\zeta}^{x}f^{\prime}(u)du\right)A_{n}(x,\zeta)d{\zeta}\nonumber\\
&&+\int_{x}^{\infty}\left(\int_{x}^{\zeta}f^{\prime}(u)du\right)A_{n}(x,\zeta)d{\zeta}.
\end{eqnarray}
Let us suppose that $$L_{1}=\int_{0}^{x}\left(\int_{t}^{x}f^{\prime}(u)du\right)A_{n}(x,\zeta)d{\zeta},$$
$$L_{2}=\int_{x}^{\infty}\left(\int_{x}^{\zeta}f^{\prime}(u)du\right)A_{n}(x,\zeta)d{\zeta}.$$
We know that $\int_{x}^{\zeta}\delta_{x}(u)d{u}=0$ and, from (\ref{eq5.2}), we have
\begin{eqnarray}\label{eq5.4}
L_{1}&=& \int_{0}^{x}\left\{\int_{\zeta}^{x}\left(\frac{1}{2}(f^{\prime}(x+)+f^{\prime}(x-))+f^{\prime}_{x}(u)+\frac{1}{2}(f^{\prime}(x+)-f^{\prime}(x-))sgn(u-x)\right)du\right\}A_{n}(x,\zeta)d{\zeta}\nonumber\\
&=&\frac{1}{2}(f^{\prime}(x+)+f^{\prime}(x-))\int_{0}^{x}(x-\zeta)A_{n}(x,\zeta)d{\zeta}+\int_{0}^{x}\left(\int_{\zeta}^{x}f^{\prime}_(x)(u)du\right)A_{n}(x,\zeta)d{\zeta}\nonumber\\
&&-\frac{1}{2}(f^{\prime}(x+)-f^{\prime}(x-))\int_{0}^{x}(x-\zeta)A_{n}(x,\zeta)d{\zeta}.
\end{eqnarray}Similarly, we find $L_{2}$
\begin{eqnarray}\label{eq5.5}
L_{2}&=& \int_{x}^{\infty}\left\{\int_{x}^{\zeta}\left(\frac{1}{2}(f^{\prime}(x+)+f^{\prime}(x-))+f^{\prime}_{x}(u)+\frac{1}{2}(f^{\prime}(x+)-f^{\prime}(x-))sgn(u-x)\right)du\right\}A_{n}(x,\zeta)d{\zeta}\nonumber\\
&=&\frac{1}{2}(f^{\prime}(x+)+f^{\prime}(x-))\int_{x}^{\infty}(\zeta-x)A_{n}(x,\zeta)d{\zeta}+\int_{x}^{\infty}\left(\int_{x}^{\zeta}f^{\prime}_{x}(u)du\right)A_{n}(x,\zeta)d{\zeta}\nonumber\\
&&+\frac{1}{2}(f^{\prime}(x+)-f^{\prime}(x-))\int_{x}^{\infty}(\zeta-x)A_{n}(x,\zeta)d{\zeta}.
\end{eqnarray}
Combining the equations (\ref{eq5.3}-\ref{eq5.4}), we get
\begin{eqnarray*}
S_{n}^{\rho}(f;x)-f(x)&=& \frac{1}{2}(f^{\prime}(x+)+f^{\prime}(x-))\int_{0}^{\infty}(\zeta-x)A_{n}(x,\zeta)d{\zeta}\\
&&+\frac{1}{2}(f^{\prime}(x+)-f^{\prime}(x-))\int_{0}^{\infty}|\zeta-x|A_{n}(x,t)d{\zeta}\\
&&-\int_{0}^{x}\left(\int_{\zeta}^{x}f^{\prime}_{x}(u)du\right)A_{n}(x,\zeta)d{\zeta}+\int_{x}^{\infty}\left(\int_{x}^{\zeta}f^{\prime}_{x}(u)du\right)A_{n}(x,\zeta)d{\zeta}.
\end{eqnarray*}Hence,
\begin{eqnarray}\label{eq5.6}
|S_{n}^{\rho}(f;x)-f(x)|&\leq& \left|\frac{f^{\prime}(x+)+f^{\prime}{x-}}{2}\right||S_{n}^{\rho}((\zeta-x);x)|\nonumber\\
&&+\left|\frac{f^{\prime}(x+)-f^{\prime}(x-)}{2}\right|S_{n}^{\rho}(|\zeta-x|;x)\nonumber\\
&&+\left|\int_{0}^{x}\left(\int_{\zeta}^{x}f^{\prime}_{x}(u)du\right)A_{n}(x,\zeta)d{\zeta}\right| \nonumber\\
&&+\left|\int_{x}^{\infty}\left(\int_{x}^{\zeta}f^{\prime}_{x}(u)du\right)A_{n}(x,\zeta)d{\zeta}\right|.
\end{eqnarray}Again, take it
$$C_{n}(f^{\prime}_{x},x)=\int_{0}^{x}\left(\int_{\zeta}^{x}f^{\prime}_{x}(u)du\right)A_{n}(x,\zeta)d{\zeta}$$ and,
$$D_{n}(f^{\prime}_{x},x)=\int_{x}^{\infty}\left(\int_{x}^{\zeta}f^{\prime}_{x}(u)du\right)A_{n}(x,\zeta)d{\zeta}.$$
Our aim is to calculate $C_{n}(f^{\prime}_{x},x),\,D_{n}(f^{\prime}_{x},x)$. From the definition $\xi_{n}(x,\zeta)$ and applying integration by parts, we get
\begin{eqnarray*}
C_{n}(f^{\prime}_{x},x)=\int_{0}^{x}\left(\int_{\zeta}^{x}f^{\prime}_{x}(u)du\right)\frac{\partial \xi_{n}(x,\zeta)}{\partial \zeta}d{\zeta}=\int_{0}^{x}f^{\prime}_{x}(\zeta)\xi_{n}(x,\zeta)d{\zeta}.
\end{eqnarray*}Thus,
\begin{eqnarray*}
|C_{n}(f^{\prime}_{x},x)|&=&\int_{0}^{x}|f^{\prime}_{x}(\zeta)|\xi_{n}(x,\zeta)d{\zeta}\\
&\leq& \int_{0}^{x-\frac{x}{sqrt{n}}}|f^{\prime}_{x}(\zeta)|\xi_{n}(x,\zeta)d{\zeta}+\int_{x-\frac{x}{\sqrt{n}}}^{x}|f^{\prime}_{x}(\zeta)|\xi_{n}(x,\zeta)d{\zeta}.
\end{eqnarray*} Since $f^{\prime}_{x}(x)=0$ and $\xi_{n}(x,\zeta)\leq 1$, we get
\begin{eqnarray*}
\int_{x-\frac{x}{\sqrt{n}}}^{x}|f^{\prime}_{x}(\zeta)|\xi_{n}(x,\zeta)d{\zeta}&=&\int_{x-\frac{x}{\sqrt{n}}}^{x}|f^{\prime}_{x}(\zeta)-f^{\prime}_{x}(x)|\xi_{n}(x,\zeta)d{\zeta}\leq \int_{x-\frac{x}{\sqrt{n}}}^{x}\left(\bigvee_{t}^{x}f^{\prime}_{x}\right)d{\zeta}\\
&\leq& \left(\bigvee_{x-\frac{x}{\sqrt{n}}}^{x}f^{\prime}_{x}\right)\int_{x-\frac{x}{\sqrt{n}}}^{x}d{\zeta}=\frac{x}{\sqrt{n}}\left(\bigvee_{x-\frac{x}{\sqrt{n}}}^{x}f^{\prime}_{x}\right).
\end{eqnarray*} From the Lemma \ref{lemma5.1} and, take $\zeta=x-\frac{x}{u}$, we have
\begin{eqnarray*}
\int_{0}^{x-\frac{x}{\sqrt{n}}}|f^{\prime}_{x}(\zeta)|\xi_{n}(x,\zeta)d{\zeta}&\leq& C_{1}|\varpi(x)|\int_{0}^{x-\frac{x}{\sqrt{n}}}\frac{|f^{\prime}_{x}(\zeta)|}{(x-\zeta)^2}d{\zeta}\\
&\leq& C_{1}|\varpi(x)|\int_{0}^{x-\frac{x}{\sqrt{n}}}\left(\bigvee_{\zeta}^{x}f^{\prime}_{x}\right)\frac{d\zeta}{(x-\zeta)^2}\\
&=&\frac{C_{1}|\varpi(x)|}{x}\int_{1}^{\sqrt{n}}\left(\bigvee_{x-\frac{x}{u}}^{x}f^{\prime}_{x}\right)\\
&\leq& \frac{C_{1}|\varpi(x)|}{x}\sum_{i=1}^{[\sqrt{n}]}\left(\bigvee_{x-\frac{x}{i}}^{x}f^{\prime}_{x}\right).
\end{eqnarray*}Finally, we get
\begin{eqnarray*}
|C_{n}(f^{\prime}_{x},x)|=\frac{C_{1}|\varpi(x)|}{x}\sum_{i=1}^{[\sqrt{n}]}\left(\bigvee_{x-\frac{x}{i}}^{x}f^{\prime}_{x}\right)+\frac{x}{\sqrt{n}}\left(\bigvee_{x-\frac{x}{i}}^{x}f^{\prime}_{x}\right).
\end{eqnarray*} Similarly integration by parts on $D_{n}(f^{\prime}_{x},x)$ and using Lemma \ref{lemma5.1}, we obtain
\begin{eqnarray*}
|D_{n}(f^{\prime}_{x},x)|&\leq& \left|\int_{x}^{2x}\left(\int_{x}^{\zeta}f^{\prime}_{x}(u)du\right)\frac{\partial}{\partial \zeta}(1-\xi_{n}(x,\zeta))d{\zeta}\right|\\
&&+\left|\int_{2x}^{\infty}\left(\int_{x}^{\zeta}f^{\prime}_{x}(u)du\right)A_{n}(x,\zeta)d{\zeta}\right|\\
&\leq& \left|\int_{x}^{2x}f^{\prime}_{x}(u)du\right||1-\xi_{n}(x,2x)|+\int_{x}^{2x}|f^{\prime}_{x}(\zeta)|(1-\xi_{n}(x,\zeta))d{\zeta}\\
&&+\left|\int_{2x}^{\infty}(f(\zeta)-f(x))A_{n}(x,\zeta)d{\zeta}\right|+|f^{\prime}(x+)|\left|\int_{2x}^{\infty}(\zeta-x)A_{n}(x,\zeta)d{\zeta}\right|.
\end{eqnarray*} Also, we have
\begin{eqnarray}\label{eq5.7}
\int_{x}^{2x}|f^{\prime}_{x}(\zeta)|(1-\xi_{n}(x,\zeta))d{\zeta}&=&\int_{x}^{x+\frac{x}{\sqrt{n}}}|f^{\prime}_{x}(\zeta)|(1-\xi_{n}(x,\zeta))d{\zeta}\nonumber\\
&&+\int_{x+\frac{x}{\sqrt{n}}}^{2x}|f^{\prime}_{x}(\zeta)|(1-\xi_{n}(x,\zeta))d{\zeta}\nonumber\\
&&J_{1}+J_{2}.
\end{eqnarray}Since $f^{\prime}_{x}(x)=0$ and, $1-\xi_{n}(x,\zeta)\leq 1$, we have
\begin{eqnarray*}
J_{1}&=&\int_{x}^{x+\frac{x}{\sqrt{n}}}|f^{\prime}_{x}(\zeta)-f^{\prime}_{x}(x)|(1-\xi_{n}(x,\zeta))d{\zeta}\\
&\leq& \int_{x}^{x+\frac{x}{\sqrt{n}}}\left(\bigvee_{x}^{x+\frac{x}{\sqrt{n}}}f^{\prime}_{x}\right)d{\zeta}\\
&=&\frac{x}{\sqrt{n}}\left(\bigvee_{x}^{x+\frac{x}{\sqrt{n}}}f^{\prime}_{x}\right).
\end{eqnarray*} From the Lemma \ref{lemma5.1} and, assuming $\zeta=x+\frac{x}{u}$, we get
\begin{eqnarray*}
J_{2}&\leq& C_{1}|\varpi(x)|\int_{x+\frac{x}{\sqrt{n}}}^{2x}\frac{1}{(\zeta-x)^2}|f^{\prime}_{x}(\zeta)-f^{\prime}_{x}(x)|d{\zeta}\\
&\leq& C_{1}|\varpi(x)|\int_{x+\frac{x}{\sqrt{n}}}^{2x}\frac{1}{(\zeta-x)^2}\left(\bigvee_{x}^{x+\frac{x}{u}}f^{\prime}_{x}\right)du\\
&=& \frac{C_{1}|\varpi(x)|}{x}\int_{1}^{\sqrt{n}}\left(\bigvee_{x}^{x+\frac{x}{u}}f^{\prime}_{x}\right)du\leq \frac{C_{1}|\varpi(x)|}{x}\sum_{i=1}^{[\sqrt{n}]}\int_{i}^{i+1}\left(\bigvee_{x}^{x+\frac{x}{u}}f^{\prime}_{x}\right)du\\
&\leq& \frac{C_{1}|\varpi(x)|}{x}\sum_{i=1}^{[\sqrt{n}]}\left(\bigvee_{x}^{x+\frac{x}{i}}f^{\prime}_{x}\right).
\end{eqnarray*}
Substitute the values of $J_{1}$ and, $J_{2}$ in (\ref{eq5.7}), we have
\begin{eqnarray*}
\int_{x}^{2x}|f^{\prime}_{x}(\zeta)|(1-\xi_{n}(x,\zeta))d{\zeta}\leq \frac{x}{\sqrt{n}}\left(\bigvee_{x}^{x+\frac{x}{\sqrt{n}}}f^{\prime}_{x}\right)+\frac{C_{1}|\varpi(x)|}{x}\sum_{i=1}^{[\sqrt{n}]}\left(\bigvee_{x}^{x+\frac{x}{u}}f^{\prime}_{x}\right).
\end{eqnarray*} Applying the Cauchy-Schwarz inequality and Lemma \ref{lemma5.1}, we get
\begin{eqnarray}\label{eq5.8}
|D_{n}(f^{\prime}_{x},x)|&\leq& M_{f}\int_{2x}^{\infty}(1+\zeta^2)A_{n}(x,\zeta)d{\zeta}+|f(x)|\int_{2x}^{\infty}A_{n}(x,\zeta)d{\zeta}\nonumber\\
&&+|f^{\prime}(x+)|\sqrt{C_{1}|\varpi(x)|}+\frac{C_{1}|\varpi(x)|}{x^2}|f(2x)-f(x)-xf^{\prime}(x+)|\nonumber\\
&&+\frac{x}{\sqrt{n}}\left(\bigvee_{x}^{x+\frac{x}{\sqrt{n}}}f^{\prime}_{x}\right)+\frac{C_{1}|\varpi(x)|}{x}\sum_{i=1}^{[\sqrt{n}]}\left(\bigvee_{x}^{x+\frac{x}{i}}f^{\prime}_{x}\right).
\end{eqnarray} Since $\zeta\leq 2(\zeta-x)$ and, $x\leq \zeta-x$ when $\zeta\geq 2x$, we have
\begin{eqnarray}\label{eq5.9}
&&M_{f}\int_{2x}^{\infty}(1+\zeta^2)A_{n}(x,\zeta)d{\zeta}+|f(x)|\int_{2x}^{\infty}A_{n}(x,\zeta)d{\zeta}\nonumber\\
&&\leq(M_{f}+|f(x)|)\int_{2x}^{\infty}A_{n}(x,\zeta)d{\zeta}+4M_{f}\int_{2x}^{\infty}(\zeta-x)^2A_{n}(x,\zeta)d{\zeta}\nonumber\\
&&\leq\frac{M_{f}+|f(x)|}{x^2}\int_{0}^{\infty}(\zeta-x)^2A_{n}(x,\zeta)d{\zeta}+4M_{f}\int_{0}^{\infty}(\zeta-x)^2 A_{n}(x,t)d{\zeta}\nonumber\\
&&\leq \left(4M_{f}+\frac{M_{f}+|f(x)|}{x^2}\right)C_{1}|\varpi(x)|.
\end{eqnarray}Using the above inequality, we have
\begin{eqnarray}\label{eq5.10}
|D_{n}(f^{\prime}_{x},x)|&\leq& \left(4M_{f}+\frac{M_{f}+|f(x)|}{x^2}\right)C_{1}|\varpi(x)|+|f^{\prime}(x+)|\sqrt{C_{1}|\varpi(x)|}\nonumber\\
&&+C_{1}\frac{1+x^2}{nx^2}|f(2x)-f(x)-xf^{\prime}(x+)|\nonumber\\
&&+\frac{x}{\sqrt{n}}\left(\bigvee_{x}^{x+\frac{x}{\sqrt{n}}}f^{\prime}_{x}\right)+\frac{C_{1}|\varpi(x)|}{x}\sum_{i=1}^{[\sqrt{n}]}\left(\bigvee_{x}^{x+\frac{x}{i}}f^{\prime}_{x}\right).
\end{eqnarray}We obtain, the required result using (\ref{eq5.6}, \ref{eq5.8}) and, (\ref{eq5.10}).
\end{proof}

\section*{Conflict of interest} The authors declare that they have no conflict of interest.


\begin{thebibliography}{1}

\bibitem {AC1996} A. Ciupa, \textit{On the approximation by Favard-Szasz type operators}. Rev. Anal. Numér. Théor. Approx. \textbf{25} 57–61, 1996.

\bibitem{AC1994} A. Ciupa, \textit{On a generalized Favard-Szasz type operator}. Research Seminar on Numerical and Statistical Calculus, Univ. Babe ¸s-Bolyai Cluj-Napoca, preprint nr. \textbf{1} 33–38, 1994.

\bibitem{BW} B. Wood, \textit{Generalized Szasz operators for the approximation in the complex domain}. SIAM J. Appl. Math. \textbf{17 (4)} 790–801, 1969.
\bibitem{DZTV} Ditzian, Z., Totik, V. \textit{Moduli of Smoothness}. Springer Series in Computational Mathematics, \textbf{9} Springer, New York, 1987.

\bibitem{DRA} Devore, R. A., Lorentz, G. G. \textit{Constructive Approximation}. Springer: Berlin/Heidelberg, Germany, 1993.

\bibitem{GIRI} Gavrea, I., Rasa, I., \textit{Remarks on some quantitative Korovkin-type results}. Rev. Anal. Nume\'{e}r, Th\'{e}or. Approx, \textbf{22} 173-176, 1993.

\bibitem{NKGVG06} Govil N. K., Gupta V., Noor M. A., \textit{Simultaneous approximation for the Phillips operators}.
Int. J. Math. Math. Sci., pp. Art. ID 49 094, \textbf{(9)}, 2006, doi: 10.1155/IJMMS/2006/49094. [Online].
Available: https://doi.org/10.1155/IJMMS/2006/49094.    

\bibitem{VGRN06}  Gupta V., Mohapatra R. N., Finta Z., \textit{A certain family of mixed summation-integral
type operators}. Math. Comput. Modelling, \textbf{42} no. 1-2, pp. 181–191, 2005, doi:
10.1016/j.mcm.2004.02.042. [Online]. Available: https://doi.org/10.1016/j.mcm.2004.02.042
    
\bibitem{HGRP} H. Gonska, R. Pˇaltˇanea, \textit{Simultaneous approximation by a class of Bernstein--Durrmeyer operators preserving linear functions}. Czech. Math. J. \textbf{60 (135)} 783-799, 2010.

\bibitem{IMEH} Ismail, M. E. H. \textit{On a generalization of Sz$\acute{a}$sz operators}. Mathematica, \textbf{39} 259-267, 1974.

\bibitem{IGVSSS} \.{I}\c{c}\"{o}z, G. Varma, S., Sucu, S. \textit{Approximation by operators including generalized Appell polynomials}. Filomat, \textbf{30} 429-440, 2016.
    
\bibitem{EIEAG} Ibikli, E, Gadjieva, EA, \textit{The order of approximation of some unbounded functions by the sequences of positive linear operators}. Turk. J. Math, \textbf{19(3)} 331-337, 1995.

\bibitem{JALD} Jakimovski, A., Leviatan, D. \textit{Generalization Sz$\acute{a}$sz operators for the approximation in the infinite interval}. Mathematica, \textbf{11} 97-103, 1969.

\bibitem{KYA} Kazmin, Y. A. \textit{On Appell polynomials}. Mat. Zametki, \textbf{6} 161-172, 1969; English translation in Math. Notes, \textbf{5} 556-562, 1969.
\bibitem{CPM77} May C. P., \textit{On Phillips operator}. J. Approximation Theory, \textbf{20 (4)} 315–332, 1977,
doi: 10.1016/0021-9045(77)90078-8. [Online]. Available: https://doi.org/10.1016/0021-9045(77)
90078-8

\bibitem{MMSRKJA} Mursaleen M., R. Shagufta, A. K. J. \textit{Approximation by Jakimovski-Leviatan-Stancu-Durrmeyer Type operators}. Filomat. \textbf{33(6)} 1517-1530, 2019.

\bibitem{CPDKND21} Prakash, C. Verma D. K., Deo N., \textit{Approximation by a new sequence of operators involving Apostol-Genocchi polynomials}. Math. Slovaca \textbf{71(5)} 1179-1188, 2021, https://doi.org/10.1515/ms-2021-0047.

\bibitem{CPNDDK21} Prakash C., Deo N., Verma D. K., \textit{Bezier variant of  Bernstein--Durrmeyer blending-type operators}. AEJM (2021) https://doi.org/10.1142/S1793557122501030.

\bibitem{CPNDDK} Prakash C., Deo N., Verma D. K., \textit{Approximation by Apostol-Genocchi sumation-integral type operators}. Miskolc Mathematical Notes, in press.
\bibitem{RP} P\u{a}lt\u{a}nea, R., \textit{Modified Szász–Mirakjan operators of integral form}, Carpathian J. Math. \textbf{24}(3) 378–385, 2008.

\bibitem{PHLP} Phillips, R. S.,\textit{An inversion formula for Laplace transforms and semi-groups of linear operators}, Annals of Math. {\bf 59}, 325--356, 1954.
    
\bibitem{MSPNSA} Sidharth M., Agrawal P. N., Araci S. \textit{Sz$\acute{a}$sz-Durrmeyer operators involving Boas-Busk polynomials of blending type}. Journal of Inequalities and Applications, \textbf{(1)} 2017 :122. doi: 10.1186/s13660-017-1396-x.

\bibitem{OZ} Sz$\acute{a}$sz, O., \textit{Generalization of S. Bernstein's polynomials to the infinite interval}. J. Res. Natl. Bur. Stand, \textbf{97} 239-245, 1950.

\bibitem{SVSS} Serhan, V., Sezgin, S., \textit{A Generalization of Sz$\acute{a}$sz Operators by Using the Appell polynomials of Class $A^{(2)}$}. Symmetry, \textbf{14(7)} 1410, 2022; https://doi.org/10.3390/sym14071410.

\bibitem {SSSV} Sezgin, S. Serhan, V., \textit{Generalization of Jakimovski-Leviatan type Szasz operators}. Applied Mathematics and Computations, \textbf{270} 977-983, 2015.

\bibitem{ZVV} Zhuk, V. V., \textit{Functions of the type Lip1 class and S. N. Bernstein's polynomials (Russian)}, Vestnik Leningrad. Univ. Mat. Mekh. Astronom. \textbf{1} 25-30, 1989.


\end{thebibliography}
\end{document}